\documentclass[12pt]{article}

\usepackage[margin=1in]{geometry}
\usepackage{amsmath,amssymb,amsthm}
\usepackage{bm}
\usepackage{booktabs}
\usepackage{float}
\usepackage{graphicx}
\usepackage[colorlinks=true,linkcolor=blue,citecolor=blue,urlcolor=blue]{hyperref}

\providecommand{\Parasplit}{}

\title{Travel-time tomography from mean field game dynamics}
\author{%
  Longqiang Xu$^{1}$, Weishi Yin$^{1,*}$, and Hongyu Liu$^{2,*}$\\[0.5em]
  \parbox{0.92\textwidth}{\centering
    $^{1}$School of Mathematics and Statistics, Changchun University of Science and Technology, No. 7089 Weixing Road, Changchun 130022, Jilin, China\\
    $^{2}$Department of Mathematics, City University of Hong Kong, Tat Chee Avenue, Kowloon, Hong Kong, China\\
    Correspondence: \href{mailto:yinweishi@cust.edu.cn}{yinweishi@cust.edu.cn}; \href{mailto:hongyliu@cityu.edu.hk}{hongyliu@cityu.edu.hk}
  }
}
\date{May 23, 2026}

\begin{document}

\maketitle

\begin{abstract}
Travel-time tomography seeks to recover a hidden environment from external measurements generated by propagation through an anomalous region. Standard formulations treat propagation as passive, so the environment influences observations mainly by bending paths or changing travel times. Many collective systems do not operate in that regime: observed arrivals are shaped by strategic motion, congestion, and environmental costs. We formulate this active setting through mean field games, in which the unknown environment enters the running cost through a spatial cost field and observations are read from the resulting population dynamics. This yields three contributions. First, it places propagation, observation generation, and inversion within one PDE-constrained model. Second, it clarifies why the inverse problem differs structurally from passive tomography: kinetic, congestion, and environmental effects are coupled endogenously and appear through space- and time-dependent local weights rather than externally chosen global coefficients. Third, it introduces a two-stage inversion pipeline that combines diffusion-based initialization with full mean field game refinement, and numerical experiments show stable recovery under noise and across multichannel scene families. Taken together, these ingredients establish a foundational framework for a class of inverse problems in which data are generated by optimizing and interacting populations rather than by passive signals. The framework identifies the forward MFG model, admissible observation channels, structural mechanisms, and computational recovery route needed to study active tomography under collective dynamics. Such problems arise, among other settings, in biological transport, vascular flow, and subsurface groundwater dynamics.
\end{abstract}


\medskip
\noindent\textbf{ORCID.} ORCID: Longqiang Xu (0009-0004-4746-7748); Weishi Yin (0000-0001-9036-5596); Hongyu Liu (0000-0002-2930-3510).

\section{Introduction}

Travel-time tomography asks how hidden environments can be recovered from external measurements generated by propagation through an anomalous region. In classical settings, propagation is passive: rays, waves, or diffusive signals traverse the unknown region, and the inverse problem is to infer hidden environmental structure from travel-time information such as boundary first-arrival data, including incomplete-data three-dimensional settings \cite{Stefanov2019,Klibanov2019}. This viewpoint has been successful across geophysics, medical imaging, and nondestructive testing, but it does not describe settings in which observed arrivals are reshaped by decision-making, congestion, or adaptive responses to the environment.

Many collective systems do operate in exactly this active way. Pedestrian and vehicle flows navigating congestion and obstacles have often been modeled through mean-field or mean-field-game formulations \cite{LachapelleWolfram2011,BurgerDiFrancescoMarkowichWolfram2014}. Related examples include tagged-agent motion, evacuation, and minimal-time multipopulation settings, where trajectories are jointly shaped by environmental costs and strategic interactions \cite{AurellDjehiche2019,SadeghiArjmandMazanti2022,YanoKuroda2023}. State constraints and nonsmooth active motion further reinforce that the data-generating process is not passive transport but adaptive collective propagation \cite{SadeghiArjmandMazanti2022b}. Related inverse questions also arise in biological transport, vascular flow, and subsurface groundwater dynamics, where hidden structures must be inferred from observed collective or distributed motion. In such systems, observations are not merely traces of passive transport.
\Parasplit They are signatures of a coupled process in which individual optimization, density redistribution, and environmental structure evolve together. The unknown environment therefore affects data not only by altering paths, but also by reshaping the collective evolution that generates the data.

Mean field games (MFGs) provide a natural mathematical language for this regime. By coupling a representative-agent optimization problem to a population evolution equation, MFGs connect microscopic control decisions to macroscopic density flows \cite{LasryLions2007,Huang2006}. Over the last two decades they have become a central framework in stochastic control, partial differential equations, and collective dynamics. Planning formulations and their variational interpretations have clarified how global trajectory coordination can be encoded at the PDE level \cite{AchdouCamilliCapuzzo2012,Porretta2013,Porretta2014}. Other developments have addressed segregation, multipopulation interactions, and velocity-alignment effects in collective motion \cite{Cirant2015,AchdouBardiCirant2017,SantambrogioShim2021}. Higher-order and state-constrained models have further expanded the framework toward controlled active propagation \cite{AchdouMannucciMarchiTchou2020,AchdouMannucciMarchiTchou2022}.

On the computational side, important work has developed numerical solvers for evolutive, congested, and planning-type MFG systems \cite{Achdou2010,AchdouCamilliCapuzzo2012,AchdouPorretta2018,EvangelistaFerreiraGomesNurbekyanVoskanyan2018}. More recent studies have pushed these ideas toward high-dimensional and learning-based settings \cite{Ruthotto2020,Lin2021,LiuJacobsLiNurbekyanOsher2021}. Related studies have also treated constrained or higher-order dynamics that are especially relevant for active propagation settings \cite{SadeghiArjmandMazanti2022,SadeghiArjmandMazanti2022b,AchdouMannucciMarchiTchou2020,AchdouMannucciMarchiTchou2022}.

Closely related inverse-problem studies have reconstructed ground metrics, interaction kernels, running costs, and interaction energies from population observations or partial boundary measurements \cite{DingLiOsherYin2022,ChowFungLiuNurbekyanOsher2023,Liu2022}. Other inverse or coefficient-identification perspectives have appeared in mean-field control, stationary decoding, and nonlinear parabolic MFG-type systems \cite{KachrooAgarwalSastry2016,LiuLo2025,LiuLoZhang2025,KlibanovLiYang2025}. Recent studies have also developed equilibrium-correction iterations for inverse MFGs, bilevel optimization formulations for inverse mean-field games, reconstruction methods for state-independent cost functions, and uniqueness and Lipschitz-stability results for retrospective analysis of MFG systems \cite{YuLiuZhao2025,YuXiaoChenLai2024,RenSoedjakWangZhai2024,KlibanovAverboukh2024}. These studies have substantially advanced MFG inverse problems, but they do not organize propagation, observation generation, and inversion from the outset as one tomography problem driven by active collective motion. The present work targets that different class of inverse problems.

\begin{figure}[H]
\centering
\includegraphics[width=0.68\linewidth]{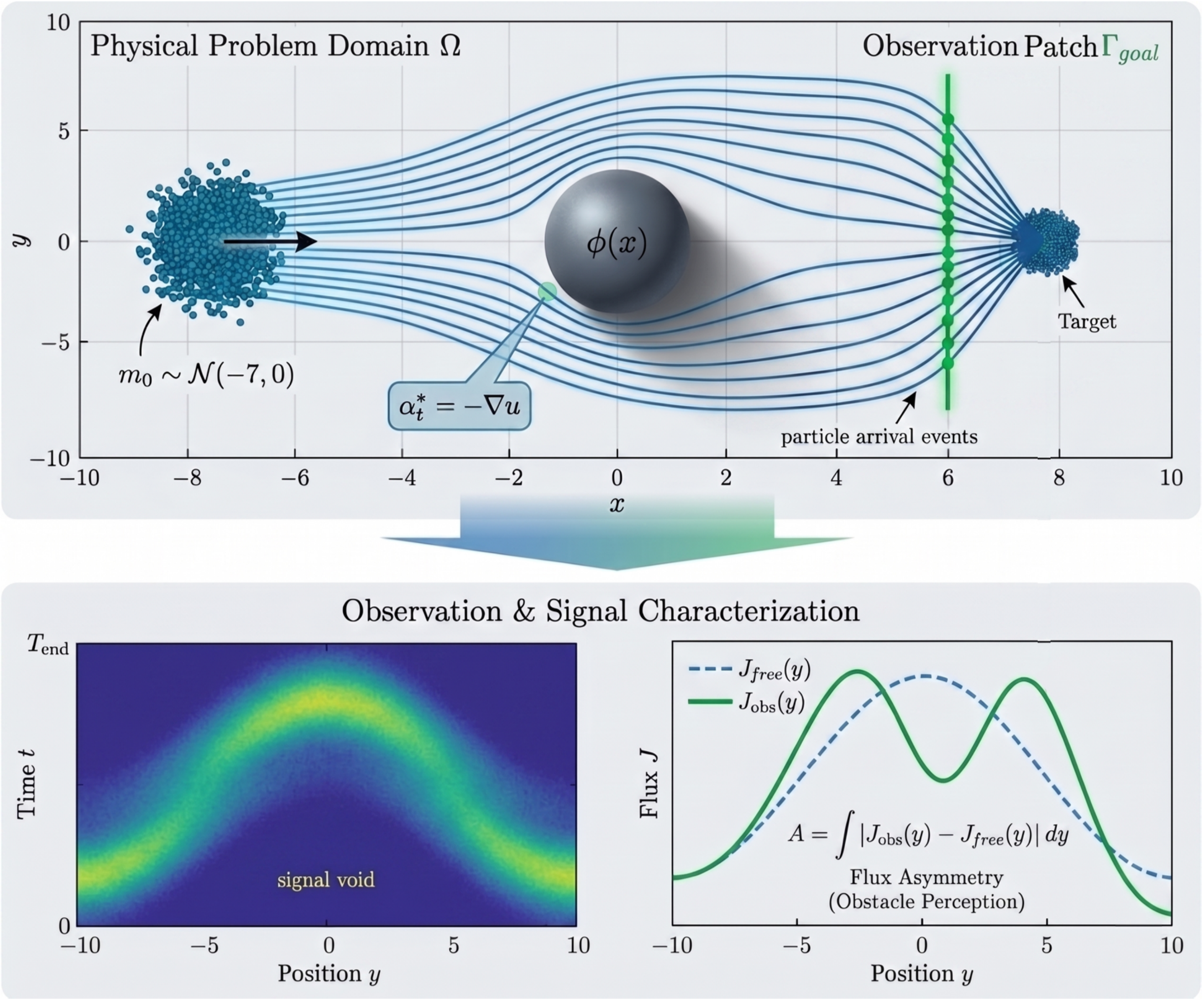}
\caption{Conceptual overview of active tomography under mean field game dynamics. The unknown environment, encoded by $\phi(x)$, alters collective propagation, which in turn generates time-resolved travel-time data, complementary terminal redistribution readouts, and derived summary statistics. Unlike passive tomography, the environment affects the data through redistribution of the population as well as through path geometry.}
\label{fig:concept}
\end{figure}

The framework studied here formulates travel-time tomography under active collective propagation. The key point is that an unknown environment enters the cost structure of the forward game through a spatial cost field and therefore influences not only optimal trajectories, but also the observation mechanism itself. Time-resolved travel-time arrival profiles, terminal redistribution data over an observation region, and low-dimensional travel-time summary statistics are treated as complementary observation channels extracted from the same forward evolution. This keeps the focus on a single tomography problem rather than recasting the manuscript as a generic multichannel fusion study. A schematic summary of this setting is shown in Fig.~\ref{fig:concept}.

This formulation yields three main contributions. First, it places propagation, observation generation, and inversion within one PDE-constrained model built from a coupled Hamilton--Jacobi--Fokker--Planck system and the associated observation dataset. Second, it clarifies why the inverse problem differs structurally from passive tomography: kinetic, congestion, and environmental effects are coupled endogenously through space- and time-dependent local weights, while stationary identifiability and long-time behavior provide structural guidance on what information can be encoded in the data. Third, it introduces a practical two-stage inversion pipeline that combines diffusion-based initialization with refinement under the full MFG model and supports stable recovery under noise and across multichannel scene families.

\section{Results}

\subsection{Active tomography and observation dataset}
We consider a bounded spatial domain $\Omega$ and an unknown environment represented by a cost field $\phi(x)$, which may encode obstacles, heterogeneous resistance, or other environmental penalties. A representative agent evolves on a finite time horizon $[0,T]$ according to a controlled diffusion,
\begin{equation}\label{eq:1}
  dX_t = \alpha_t\,dt + \sqrt{2\nu}\,dW_t,
\end{equation}
with initial distribution $m_0$. The control $\alpha_t$ is chosen to minimize a cost functional of the form
\begin{equation}\label{eq:2}
  J(\alpha)=\mathbb{E}\Big[\int_0^T \big(\tfrac{1}{2}|\alpha_t|^2 + f(m_t) + \phi(X_t)\big)\,dt + g(X_T)\Big].
\end{equation}

The three running-cost components have distinct roles: the kinetic term $|\alpha_t|^2$ penalizes inefficient motion, the congestion term $f(m_t)$ captures interactions through the population density, and the environment term $\phi(X_t)$ represents the hidden environment to be recovered. In the mean field limit, the forward propagation problem is described by a coupled Hamilton--Jacobi--Fokker--Planck system for the value function $u(x,t)$ and density $m(x,t)$,
\begin{equation}\label{eq:3}
\begin{cases}
-\partial_t u - \nu \Delta u + \tfrac{1}{2}|\nabla u|^2 = f(m) + \phi(x), \\
\partial_t m - \nu \Delta m - \nabla \cdot (m \nabla u) = 0,
\end{cases}
\end{equation}
with initial-terminal conditions supplied by $m(x,0)=m_0(x)$ and $u(x, T)=g(x)$. Since the motion is confined within the domain $\Omega$, we can impose periodic boundary conditions on $\partial\Omega$. This system is the forward model from which observations are generated. The inverse problem is to recover the hidden environment, encoded by $\phi$, from data extracted from the corresponding solution pair $(u,m)$.

For a typical and representative tomographic scenario, we take $\phi$ to be a nonnegative obstacle potential of the form
\begin{equation}\label{eq:4}
  \phi(x)=C\,\chi_{\mathcal O}(x)+\phi_{bg}(x),
\end{equation}
where $\mathcal O$ is the support of the hidden object, $C$ is the obstacle strength, and $\phi_{bg}$ is a known background field. This parameterization makes explicit that the inverse problem is not simply to estimate an abstract coefficient, but to recover a geometrically meaningful environment from population-level observations.

In classical tomography, one typically starts with a prescribed observation operator acting on a passive propagation model. Here the observation mechanism is itself part of the MFG formulation. Once the forward system is solved for a given environment field $\phi$, one may extract several observation channels from the same evolution: time-resolved arrival profiles on the target boundary, terminal redistribution patterns over an observation region, or low-dimensional statistics derived from the time-dependent arrival history.

This distinction is important for this work's central claim. The forward problem remains one coupled propagation system, while multiple observations correspond only to different travel-time readouts of that system. In particular, multi-observation settings do not introduce independent physical mechanisms; they organize different measurements of the same collective dynamics.

We therefore write the forward map schematically as
\begin{equation}\label{eq:5}
  \phi \mapsto (u,m) \mapsto \mathcal{M}(\phi),
\end{equation}
where $\mathcal{M}$ may be scalar- or vector-valued depending on the chosen observation channels. The central tomography problem is then to infer $\phi$ from noisy realizations of $\mathcal{M}(\phi)$.

The basic single-channel travel-time observation is the time-dependent arrival trace on an  observation boundary patch $\Gamma_{goal}$ (see Fig.~\ref{fig:concept} for a schematic illustration),
\begin{equation}\label{eq:6}
  \mathcal F(\phi):=\mathrm{Trace}_{\Gamma_{goal}}\circ \mathrm{Proj}_m \circ \mathcal S(\phi),
\end{equation}
where $\mathcal S$ denotes the forward MFG solver and $\mathrm{Proj}_m$ extracts the density component. Indeed, it can be explicitly cast as:
\begin{equation}\label{eq:7}
\mathcal{F}(\phi)=m(x, t)\big|_{(x, t)\in\Gamma_{goal}\times (0, T)}. 
\end{equation}
For the tomography interpretation emphasized here, this boundary trace is not treated merely as a boundary density measurement. Rather, it is the collective travel-time signal in the present active travel-time tomography setting: its onset, temporal concentration, and accumulated mass encode how long the controlled population takes to reach the target region through the hidden environment. Thus \eqref{eq:6} supplies the travel-time observation used in the single-channel setting, with the data represented as a full time-resolved arrival profile rather than a single passive first-arrival scalar. As illustrated in Fig.~\ref{fig:concept}, the arrival events on the observation boundary patch generate a time-resolved signal. 

If the observation dataset is collected from a single initial population distribution \( m_0 \), it is referred to as a single measurement; otherwise, multiple measurements are performed. This paper focuses primarily on the single-measurement case. Such tomographic problems arise in various settings, including biological transport, vascular flow, and subsurface groundwater dynamics.

\subsection{Endogenous structure and identifiability}
One of the key structural features of the formulation is that the different physical effects are coupled endogenously rather than by hand-tuned global weights. The forward MFG system can be related to a constrained variational problem in which density evolution and velocity fields are optimized jointly. In this interpretation, kinetic costs, congestion penalties, and environmental costs are coordinated through a common optimality structure.

This matters for the inverse problem because it clarifies what information is genuinely encoded in the data. A local change in the recovered environment is not interpreted solely as a geometric deflection of trajectories. It also modifies how kinetic efficiency and congestion avoidance are balanced during the collective evolution. The resulting observations therefore contain signatures of the environment both directly and through the redistribution of the population. This endogenous coupling is one of the main reasons the MFG setting differs from classical passive tomography and from inverse formulations built from externally weighted surrogate models.

This idea can be made more concrete through local effective weights. At an equilibrium $(u^\ast,m^\ast)$, the kinetic, congestion, and environmental contributions define pointwise ratios such as
\begin{equation}\label{eq:8}
  w_{kin}(x,t)=\frac{\tfrac12 |\nabla u^\ast|^2}{\tfrac12 |\nabla u^\ast|^2+f(m^\ast)+\phi(x)},
\end{equation}
\begin{equation}\label{eq:9}
  w_{crowd}(x,t)=\frac{f(m^\ast)}{\tfrac12 |\nabla u^\ast|^2+f(m^\ast)+\phi(x)},
\end{equation}
with an analogous expression for the environment term. These weights are not externally chosen hyperparameters. They are induced by the equilibrium itself, vary over space and time, and help explain why the data encode a redistribution process rather than a single fixed tradeoff between competing effects.

\begin{figure}[H]
\centering
\includegraphics[width=0.85\textwidth]{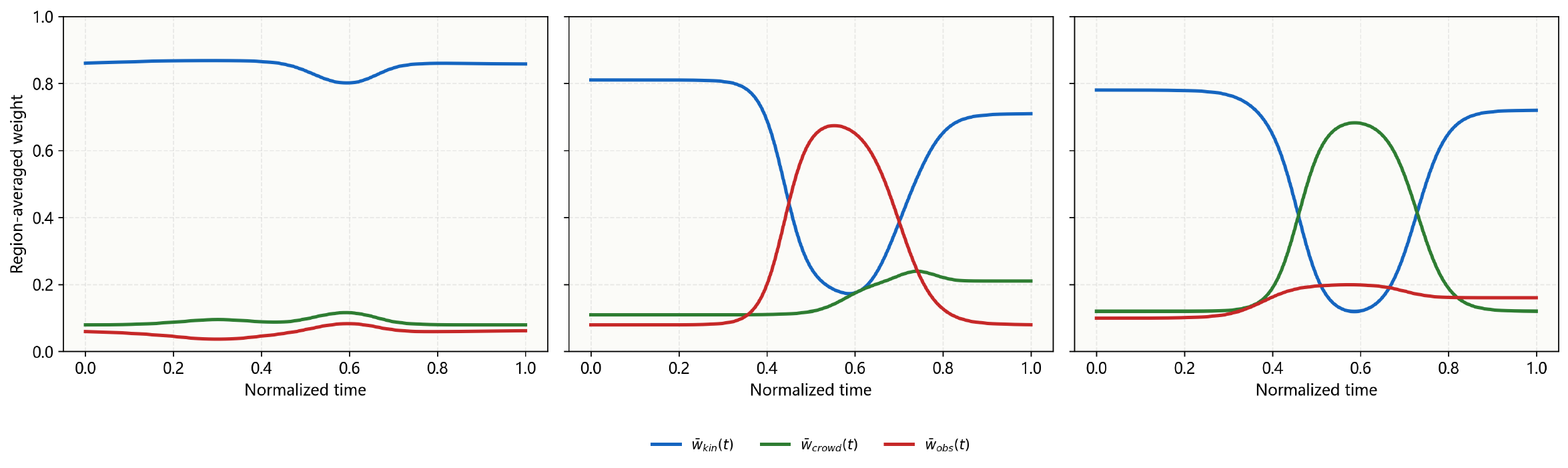}
\caption{Mechanism-oriented view of endogenous local weights. The three panels show region-averaged temporal weights in an open corridor, near an obstacle, and in a bottleneck regime. Efficient motion dominates in open regions, obstacle cost rises sharply near barriers, and congestion becomes most prominent under compression.}
\label{fig:weights}
\end{figure}

Two-dimensional mechanism experiments sharpen this interpretation. In open regions, $w_{kin}$ dominates along the main propagation corridor, reflecting a preference for efficient travel. Near strong obstacles, the environment-related weight rises locally as the population approaches the barrier, consistent with rapid lateral redistribution during obstacle avoidance. In bottleneck regions formed by paired obstacles, the congestion-related weight becomes dominant in the densest part of the flow, indicating that crowd interaction is activated most strongly when local compression is highest. These observations support the claim that the MFG formulation does not impose one global balance across the whole domain; instead, it generates a spatially and temporally varying response profile. The corresponding temporal weight profiles are illustrated in Fig.~\ref{fig:weights}.

The dynamic inverse problem \eqref{eq:7} is challenging because the observations are generated by a nonlinear forward-backward system and are only available on a boundary patch over a finite time window. A primary theoretical concern is unique identifiability, which asks for sufficient conditions guaranteeing the one-to-one correspondence:
\begin{equation}\label{eq:10}
  \mathcal F(\phi_1)=\mathcal F(\phi_2) \iff \phi_1=\phi_2. 
\end{equation}
When multiple measurements are available, existing studies \cite{LiuLo2025,LiuLoZhang2025} provide a theoretical roadmap to establish unique identifiability in generic scenarios. However, using only a single measurement introduces significant difficulties. This can be seen intuitively from a cardinality argument: the observational data are collected on a boundary patch, whereas the unknown is a potential function supported inside a bulk body. Nevertheless, if the goal is restricted to recovering only the geometric shape and location of the hidden inclusion $\mathcal{O}$, the problem becomes highly feasible. This is especially true when the obstacle strength $C$ in \eqref{eq:4} is sufficiently large, in which case the motion does not enter $\mathcal{O}$. We will report these theoretical findings in a forthcoming work. In the present paper, we focus on laying the foundational framework for this active and collective travel-time tomography derived from MFG dynamics.

\subsection{Two-stage inversion}
The computational strategy is organized as a two-stage inversion method tailored to the structure and computational difficulty of active tomography under mean field game dynamics. This design is not an arbitrary algorithmic decomposition. A direct one-stage inversion under the full active MFG model is strongly nonlinear, tightly coupled, and highly sensitive to initialization. The two-stage structure is introduced to separate coarse environmental recovery from full-model refinement while preserving the active propagation logic of the problem.

In the first stage, a reduced inverse model provides a low-frequency initialization for the hidden environment field. This stage is motivated by the characteristic terminal-approach behavior of agents in the MFG system. Under suitable continuity assumptions, when agents approach the target region their control is dominated primarily by target-entry and congestion effects. Along directions parallel to the target boundary, the effective directed-control component becomes much weaker, so the local dynamics are well approximated by Brownian motion coupled with crowding. This reduced terminal-approach regime motivates a diffusion-based coarse inversion stage: it does not replace the full MFG model, but extracts large-scale geometric information at lower computational cost and reduces the instability of blind initialization in a strongly nonlinear inverse problem.

In the second stage, the initialization is refined under the full MFG-based inverse formulation. This is the stage at which the endogenous coupling among motion, congestion, and environmental structure is restored and becomes fully relevant to reconstruction. Stage I therefore captures the coarse environmental signature under a reduced terminal-approach surrogate, whereas Stage II refines that estimate under the complete active propagation model. In low-dimensional settings, this refinement is carried out directly under the coupled Hamilton--Jacobi--Fokker--Planck constraint. In higher-dimensional settings, the same refinement role may be realized through scalable surrogate inversion modules, provided that they are still guided by the full MFG forward structure and observations generated by the same active propagation model.

In the notation adopted here, a canonical single-channel realization of the second stage minimizes an objective of the form
\begin{equation}\label{eq:11}
  \mathcal L(\phi)=\frac12 \int_0^T \int_{\Gamma_{\text{goal}}} |m(x,t;\phi)-y_{\text{obs}}(x,t)|^2\,dx\,dt+\frac{\lambda}{2}\|\phi\|_{L^2}^2,
\end{equation}
subject to the full mean-field game (MFG) system. When multichannel observations are available, this mismatch functional can be modified by a straightforward vector-valued extension, replacing the scalar field \(m\) and scalar observation \(y_{\text{obs}}\) with their vector counterparts. Moreover, in multichannel implementations, the same two-stage logic can be preserved while replacing direct full-model refinement with scalable surrogate realizations \cite{Ruthotto2020,Lin2021,LiuJacobsLiNurbekyanOsher2021,Lauriere2023}. It is remarked that as noted earlier, we will be primarily concerned with recovering only the support \(\mathcal{O}\) of the obstacle, rather than the full potential function, especially for the single-channel realization considered in the remainder of this paper.   

This cascade strategy is especially natural in the present context because it respects the separation between coarse transport information and fine strategic effects. It also provides a practical bridge between mathematically structured PDE models and higher-dimensional learning-based implementations used in numerical experiments.

\begin{figure}[H]
\centering
\includegraphics[width=0.98\textwidth]{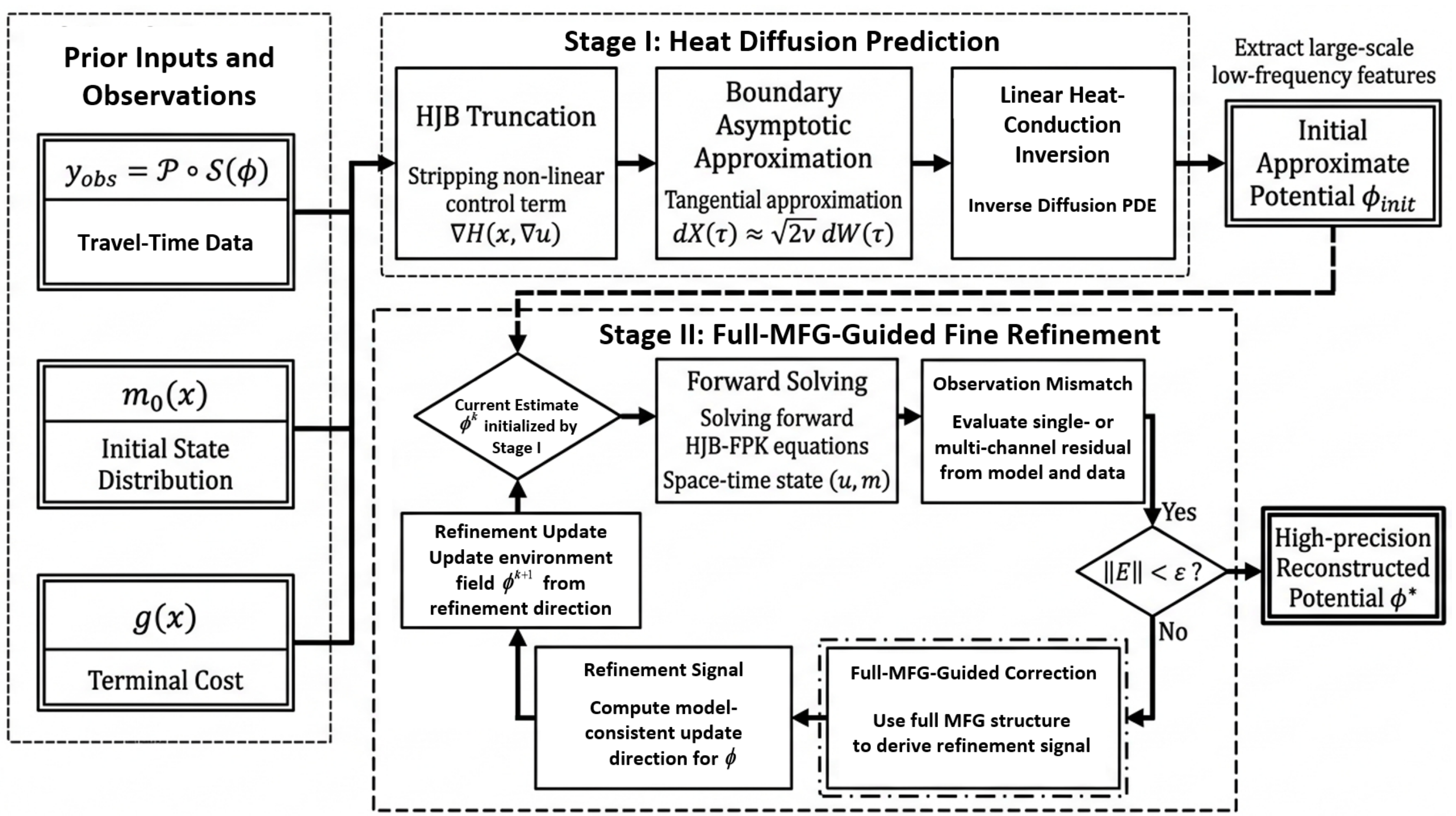}
\caption{Two-stage inversion pipeline. Stage I extracts large-scale geometric information through a diffusion-based approximation, producing a stable initialization of the hidden environment field. Stage II refines that estimate under the full Hamilton--Jacobi--Fokker--Planck constraint using forward solves, observation mismatch evaluation, model-consistent correction signals, and iterative parameter updates.}
\label{fig:pipeline}
\end{figure}

\subsection{Numerical reconstruction and cross-scene stability}
The numerical studies are organized in three roles. First, low-dimensional experiments provide a mechanistic validation of the proposed two-stage inversion method. Second, a higher-dimensional realization tests whether the same coarse-to-fine logic can be implemented in a scalable form. Third, multichannel experiments examine how the observation-side extension of the framework behaves across heterogeneous scene families.

The low-dimensional experiments use two-dimensional synthetic environments with shared forward setting $\Omega=[-2,6]\times[-2,6]$, diffusion coefficient $\nu=0.15$, horizon $T=2.0$, a $35\times 35$ spatial grid, 50 time steps, and a smooth circular obstacle centered near $(1.5,1.5)$. These tests serve as a mechanistic validation of the proposed inversion structure: they examine whether the reduced first stage supplies a usable coarse estimate under the basic dynamic operator $\mathcal F_T$, whether the full-model second stage can refine that estimate reliably, and how the reconstruction changes under finite observation windows and moderate noise. Noise-stress tests use both additive and multiplicative perturbations. Additive noise levels of $0\%$, $2\%$, and $5\%$ are combined with multiplicative noise levels of $0\%$, $10\%$, and $20\%$. Even in the strongest displayed disturbance settings, reconstructions show gradual drift rather than catastrophic topological failure. Sensitivity studies also suggest that performance stabilizes once the observation horizon is sufficiently long, consistent with the idea that longer trajectories expose more of the active propagation structure.

To assess scalability, we then consider a three-dimensional realization of the same single-observation logic. In this setting, the second-stage refinement is implemented by a lightweight surrogate module trained on MFG-generated observations. The resulting reconstruction recovers obstacle center coordinates with an average $R^2$ of 0.8824, mean absolute errors on the order of $0.10$ to $0.14$, and a radius mean absolute error of 0.0459. The implementation remains relatively compact, with about $1.44\times 10^5$ trainable parameters ($\sim 0.55$ MB) and reported inference latency on the order of $0.2$ to $0.46$ ms per reconstruction in the tested setup. Blind-test statistics further indicate a median center error of 0.1707 and a 90th-percentile error of 0.5645, while training-curve comparisons across five random seeds point to similarly stable convergence behavior. Taken together, these results suggest that the proposed coarse-to-fine inversion logic can be realized in a fast and compact high-dimensional setting without abandoning the underlying active-propagation interpretation.

We next turn to the multichannel extension of the travel-time tomography framework. Relative to the canonical single-channel setting above, these experiments combine three travel-time-oriented channels generated by the same dynamic evolution: the time-dependent boundary arrival profile, the terminal redistribution pattern over an interior observation region, and a low-dimensional summary vector involving first-arrival time, asymmetry, and accumulated mass. The goal is not to claim that adding channels automatically improves reconstruction in every setting. Rather, it is to test whether the MFG-based propagation and observation design can support multichannel travel-time inversion without scene-specific manual weight tuning.

These channels can be written as
\begin{equation}\label{eq:12}
  \mathcal F_1(\phi):=m_\phi|_{\Gamma_{goal}\times[0,T]},
  \qquad
  \mathcal F_2(\phi):=m_\phi(\cdot,T)|_{\Omega_{obs}},
\end{equation}
and
\begin{equation}\label{eq:13}
  \mathcal F_3(\phi):=\big[t_{\mathrm{first}}(\phi),\,\eta_{\mathrm{asym}}(\phi),\,M_{\mathrm{acc}}(\phi)\big]^\top,
\end{equation}
where $\Omega_{obs}\subset\Omega$ denotes an interior observation region near the target boundary, for example a thin terminal neighborhood such as $\Omega_{obs}=\{x\in\Omega:\mathrm{dist}(x,\Gamma_{goal})\leq \delta\}$, and the three entries of $\mathcal F_3$ represent first-arrival time, arrival asymmetry, and accumulated arrival mass. Thus $\mathcal F_1$ is the same time-resolved travel-time arrival-profile operator introduced in \eqref{eq:6}, now used as the first channel in the multichannel setting. The second channel records the terminal redistribution pattern associated with the same travel-time process, rather than another boundary time trace. The third is a derived low-dimensional travel-time summary of the arrival history, included to test whether coarse arrival statistics improve robustness after channel normalization.

The reported multichannel results show a clear cross-scene pattern. Relative to the single-observation MFG setting, the three-channel MFG setting reduces average center error from 0.254 to 0.195 and improves average intersection-over-union from 0.829 to 0.890. At the same time, the cross-scene standard deviation stays low at 0.043, compared with 0.122 for a traditional fixed-weight three-channel baseline. Obstacle-dominated settings, congestion-dominated settings, and mixed-conflict settings are all considered, and across all three families the MFG-based multichannel inversion shows a more monotone improvement as channels are added. For example, the reported average reconstruction error in obstacle-dominated settings decreases from 0.241 to 0.207 to 0.184 as one moves from one to two to three channels; analogous monotone decreases appear in congestion-dominated settings (0.271 to 0.233 to 0.208) and mixed-conflict settings (0.323 to 0.281 to 0.233). By contrast, the traditional fixed-weight baseline may improve from one to two channels but then degrade after a third channel is introduced.

\begin{figure}[H]
\centering
\includegraphics[width=0.48\textwidth]{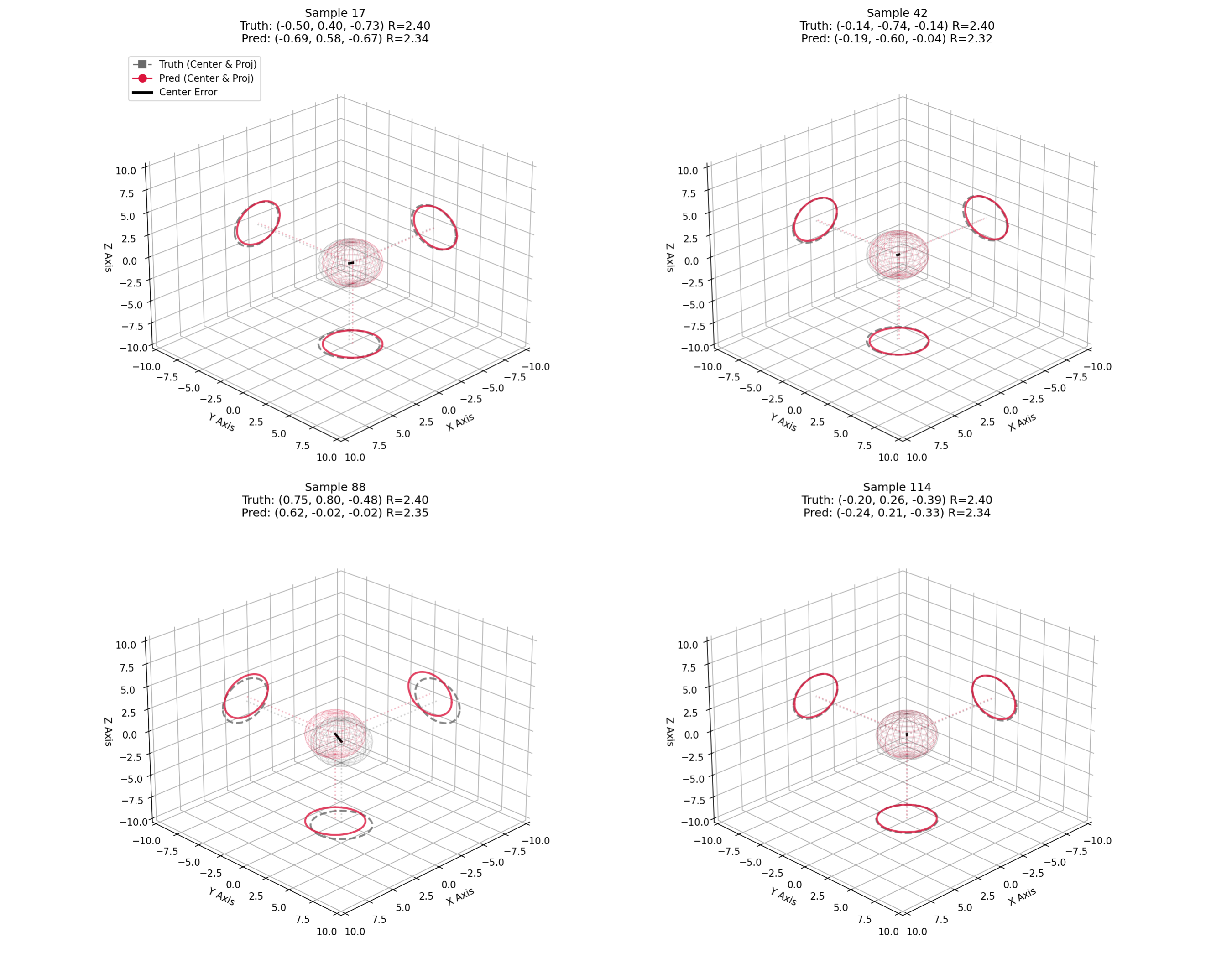}\hfill
\includegraphics[width=0.48\textwidth]{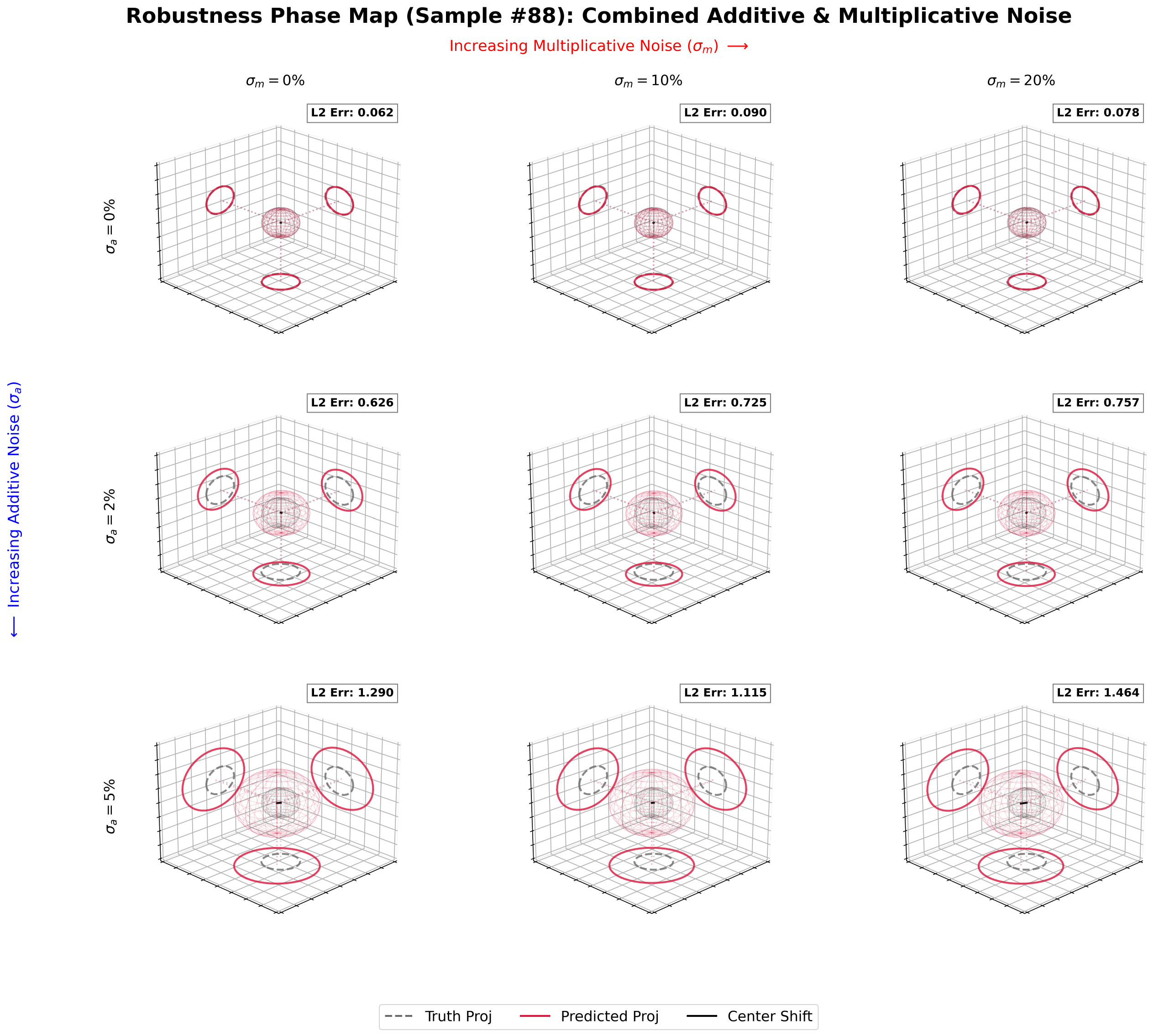}

\vspace{0.75em}
\includegraphics[width=0.72\textwidth]{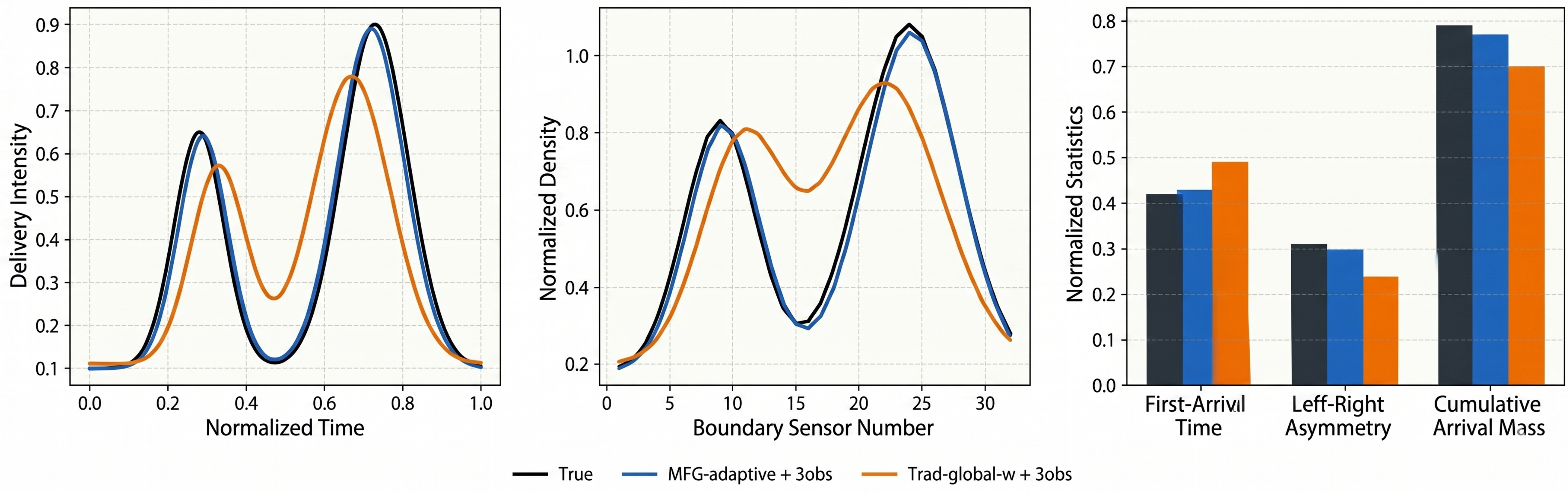}
\caption{Representative numerical evidence. Top left, sample reconstructions in the three-dimensional single-observation setting show agreement between predicted and ground-truth obstacle geometry. Top right, robustness maps indicate graceful degradation under combined additive and multiplicative noise. Bottom, multichannel comparisons show that MFG-based observation design uses three channels more stably than a traditional fixed-weight baseline across representative scene families.}
\label{fig:numerics}
\end{figure}

The contrast between the two traditional baselines is similarly informative. The fixed-weight baseline, denoted Trad-global-$\omega$, uses one observation-weight vector for all scenes. The oracle baseline, Trad-oracle-$\omega$, searches for a better weight vector separately for each scene. The fact that the oracle baseline can approach the MFG result in some cases suggests that the issue is not lack of information in the channels themselves. Rather, the issue is how robustly that information can be combined without repeated manual retuning. This is why the cross-scene stability numbers matter as much as the mean reconstruction error, and why the multichannel MFG design is best interpreted as offering a better balance between reconstruction quality and stability under scene variation rather than as an unconditional dominance claim. Representative single-observation recovery, robustness, and multichannel comparisons are summarized in Fig.~\ref{fig:numerics}.

\section{Discussion}

The main message is that tomography changes character when the data are generated by adaptive collective propagation rather than passive transport. Mean field games provide a principled model for this regime because they combine environmental costs, strategic motion, and density evolution in one forward description. The resulting inverse problem is therefore not simply parameter estimation for a fixed PDE; it is recovery of hidden environments whose observational signatures are shaped by collective behavior.

Seen through the three contributions emphasized in the abstract, the framework makes three connected points. First, propagation, observation generation, and inversion belong to one PDE-constrained model rather than to separate modeling stages. Second, the inverse problem differs structurally from passive tomography because kinetic, congestion, and environmental effects are coupled endogenously through space- and time-dependent local weights instead of being combined through externally chosen global coefficients. Third, the two-stage inversion pipeline shows that this structural viewpoint can be turned into a practical computational strategy, with numerical evidence of stable recovery under noise and across multichannel scene families.

The framework also sharpens the roles of the main ingredients. Multi-observation data are interpreted as multiple travel-time-oriented channels generated by a common propagation process rather than as unrelated information sources. Variational structure explains why kinetic, congestion, and environmental effects are coupled without externally chosen global weights. Long-time behavior clarifies why stationary inverse results remain relevant even when the central problem is dynamic. Together these points explain why classical intuition from passive travel-time tomography cannot simply be transferred to active propagation settings.

Several limitations are equally important. The strongest identifiability statements currently come from stationary or closely related settings \cite{LiuLoZhang2025} rather than from the full dynamic inverse problem. Large-scale numerical reconstruction also remains computationally demanding because each inverse step depends on solving a nonlinear coupled PDE system, even though recent numerical and learning-based MFG solvers have begun to improve scalability \cite{Ruthotto2020,Lin2021,LiuJacobsLiNurbekyanOsher2021}. These limitations mark the main directions in which the theory and computation need to mature.

The broader implication is that inverse problems can change character when the objects that generate the data are themselves optimizing, interacting, and redistributing collectively. Similar issues may arise in biological transport, vascular flow, subsurface groundwater dynamics, and other systems where hidden environments must be inferred from collective or distributed motion. In that sense, MFG-based tomography can be viewed as a model case for inverse problems under collective dynamics.

\section{Materials and Methods}

\subsection{Model assumptions}
The forward model assumes a bounded domain $\Omega$, an initial density $m_0$, a terminal cost $g$, a diffusion coefficient $\nu > 0$, and a congestion term $f$ satisfying standard monotonicity and regularity conditions. The unknown environment is encoded by a nonnegative cost field $\phi(x)$, interpreted as an obstacle or heterogeneous environmental potential. Regularity assumptions are chosen so that the coupled Hamilton--Jacobi--Fokker--Planck system is well posed in the regime under study.

\subsection{Forward and observation operators}
For each admissible environment field $\phi$, the forward operator maps $\phi$ to the MFG solution pair $(u,m)$. Observation operators then act on $(u,m)$ to produce travel-time data channels used in the inverse problem. Typical examples include time-resolved travel-time arrival profiles on a target portion of the boundary, terminal redistribution patterns over an observation region, and low-dimensional summaries of the time-dependent travel-time signal. The multi-observation setting corresponds to vector-valued observation operators built from the same forward evolution.

In the concrete three-channel setting, the observation tuple consists of a time-dependent boundary travel-time arrival profile, terminal redistribution data over an interior observation region, and a low-dimensional statistics vector encoding first-arrival time, asymmetry, and accumulated arrival mass. Here the interior observation region may be chosen as a thin terminal neighborhood of the target boundary. These channels are standardized before joint inversion so that scale differences do not dominate the mismatch functional.

\subsection{Variational and structural analysis}
The coupled PDE system can be related to a constrained variational problem for density and velocity fields,
\begin{equation}\label{eq:14}
  \min_{m,v}\int_0^T \int_\Omega \Big[\tfrac12 m|v|^2+F(m)+\phi(x)m\Big]\,dx\,dt+\int_\Omega G(m(T))\,dx,
\end{equation}
subject to the continuity equation with diffusion. This formulation is used to interpret kinetic, congestion, and environmental effects as endogenously coupled components of one optimization structure. Stationary inverse results and long-time asymptotics are then used as structural reference points for understanding the dynamic recovery problem.

\subsection{Inversion pipeline}
The inversion pipeline follows the two-stage design described in Results. The first stage is motivated by a reduced terminal-approach regime of the MFG dynamics. Under suitable continuity assumptions, when agents approach the target region their control is dominated mainly by target-entry and crowding effects, while the directed component along tangential directions becomes weaker. This motivates a reduced diffusion-based surrogate in which Brownian motion and crowding provide the leading-order coarse inversion signal. The resulting first-stage estimate is therefore intended to recover low-frequency environmental structure rather than to replace the full inverse problem.

The second stage restores the full active propagation model and refines the coarse estimate under the complete MFG-based inverse formulation. In low-dimensional settings this refinement is carried out directly under the coupled HJB--FPK constraint, whereas in higher-dimensional settings the same role may be realized through scalable surrogate modules trained on MFG-generated observations. In both cases, the methodological distinction is the same: Stage I supplies a reduced-model coarse estimate, and Stage II performs full-model refinement.

For multichannel experiments, a standardized joint loss of the form
\begin{equation}\label{eq:15}
  \mathcal L_{\mathrm{MFG}}(\phi)=\|\widehat{\mathcal Y}(\phi)-\widehat{\mathcal Y}^{obs}\|_2^2+\alpha\|\phi\|_{L^2}^2,
\end{equation}
is used, where $\widehat{\mathcal Y}$ stacks the normalized observation channels. The traditional comparison baseline instead introduces manual channel weights $\omega_1,\omega_2,\omega_3$ in the loss. This distinction is methodologically important: both approaches may preprocess the channels, but only the traditional baseline requires explicit scene-sensitive observation weighting inside the inverse objective.

\subsection{Numerical settings}
The baseline two-dimensional experiments use synthetic obstacle-reconstruction settings in which the initial density is Gaussian, the target region lies near the opposite corner of the computational domain, and observations are generated from noise-free forward solutions before additive or multiplicative perturbations are injected. The shared baseline parameters are $\nu=0.15$, $T=2.0$, a $35\times 35$ spatial grid, 50 time steps, $\lambda=1.0$, and learning rate $\eta=0.3$. Shared discretization and stopping criteria are used across methods so that differences in reconstruction quality are attributable to inversion strategy rather than to changes in the forward model.

The higher-dimensional extension uses a learning-based second-stage implementation calibrated on synthetic MFG trajectories and evaluated through obstacle-center and radius reconstruction metrics. For the higher-dimensional forward data-generation stage, synthetic MFG trajectories were produced using code adapted in part from Deep Generalized Schr\"odinger Bridge \cite{LiuChenSoTheodorou2022}. For the three-dimensional data-generation stage, the reported calibration sweep spans stage count $\texttt{num\_stage}\in\{10,20,30,40,50\}$ and per-stage iteration count $\texttt{num\_itr}\in\{50,100,150,200,250\}$, selecting $\texttt{num\_stage}=50$ and $\texttt{num\_itr}=250$ for the later experiments. The compact regression backbone used in this setting contains approximately 144,132 trainable parameters.

For multichannel experiments, three observation channels are standardized before joint inversion. The baseline comparisons then distinguish an MFG loss using normalized stacked channels from traditional weighted losses that introduce one fixed observation-weight vector for all scenes or one oracle vector retuned separately for each scene. This setup is designed specifically to test whether multichannel gains persist when scene structure changes.


\section*{Acknowledgments}
The research was supported by the Hong Kong RGC General Research Funds (projects 11311122, 12301420 and 11300821).

\end{document}